\documentstyle{amsppt}
\pagewidth{5.5in}
\pageheight{8in}
\NoRunningHeads
\magnification = 1200

\topmatter
\title HAMILTONIAN-MINIMAL SUBMANIFOLDS IN KAEHLER MANIFOLDS WITH SYMMETRIES\endtitle
\author Yuxin Dong \endauthor
\thanks {Supported by NSFC (No.10131020)}
\endthanks
\abstract
By making use of the symplectic reduction and the cohomogeneity method, we give a general method for constructing Hamiltonian minimal submanifolds in Kaehler manifolds with symmetries. As applications, we construct infinitely many nontrivial complete Hamiltonian minimal submanifolds in $CP^n$ and $C^n$.
\endabstract
\subjclass{53C42}, {53D12}
\endsubjclass
\endtopmatter
\document
\heading{\bf 1. Introduction}
\endheading
\vskip 0.3 true cm
Let $(M^{2m},\omega )$ be a symplectic manifold with a Riemannian metric $g$
and let $L$ be a Lagrangian submanifold in $M$. A normal vector field $V$
along $L$ is called a Hamiltonian variation if the one form $\alpha
_V:=\omega (V,\cdot )$ is exact. According to [O1,2], the Lagrangian
submanifold $L$ is called Hamiltonian minimal if it is a critical point of
the volume functional with respect to all Hamiltonian variations along $L$.
In particular, this makes sense if $M$ is a Kaehler manifold. A Hamiltonian
minimal submanifold will be simply called $H-$minimal.

\proclaim{Proposition 1.1}
([O2]) Let $(M,\omega ,g)$ be a Kaehler manifold. A Lagrangian submanifold $%
L\subset M$ is $H-$minimal if and only if its mean curvature vector $H$
satisfies 
$$
\delta \alpha _H=0  \tag{1}
$$
on $L$, where $\delta $ is the Hodge-dual operator of $d$ on $L$.
\endproclaim

$H-$minimal Lagrangian submanifolds offer a nice generalization of the
minimal submanifold theory. It was Oh who first investigated these
submanifolds (see [O1-2]). One motivation to study them is its similarity to
some models in incompressible elasticity([Wo], [HR1]). In [O2], the author
comments that $H-$minimal Lagrangian submanifolds seem to exist more often
than minimal Lagrangian submanifolds do. In [CU], Castro and Urbano
constructed some exotic Hamiltonian tori in $C^2$. Afterwards, Helein and
Romon constructed $H$-minimal surfaces via integrable system method
([HR1,2]). Besides these explicit instances, Schoen and Wolfson [SW]
established some important existence and regularity results for
two-dimensional $H-$minimal surfaces. However, only a few non-trivial
examples of $H-$minimal Lagrangian submanifolds of higher dimensions have
been known so far.

The aim of this paper is to give some constructions of $H-$minimal
Lagrangian submanifolds of higher dimensions. Note that the equation (1) is
a third order P.D.E., which is more complicated than the minimal submanifold
equation. Even for the usual minimal submanifold, the existence is a
difficult area of study, due to the nonlinearity of the equation. Recently,
the symmetry reduction method leads to some important progress in explicit
construction of special Lagrangian submanifolds by several authors (see
[J1,2] and the references contained therein). In this paper, we will solve
(1) by the same trick. Let $G$ be a compact connected Lie group of
holomorphic isometries of a Kaehler manifold $M$ and let $\mu $ be the
moment map of the $G-$action. First, we show that a $G-$invariant Lagrangian
submanifold is $H-$minimal if and only if it is stationary with respect to
any $G-$invariant Hamiltonian variation. From [J2], we know that a $G-$
invariant Lagrangian submanifold is contained in a level set of $\mu $. The
well-known Noether theorem tells us that the moment map $\mu $ is a
conserved quantity for every $G-$invariant Hamiltonian deformation. This
allows us to restrict the variational problem in a level set of $\mu $. By
combining the symplectic reduction and the cohomogeneity method developed in
[HsLa], we can reduce the equation (1) to a P.D.E. on the symplectic
quotient with the Hsiang-Lawson metric. We have a very nice correspondence
between the $G-$invariant $H-$minimal Lagrangian submanifolds in $M$ and the 
$H-$minimal Lagrangian submanifolds in the quotient space (see Corollary
2.8 and Theorem 2.9). The reduction procedure simplifies the original
equation greatly. Actually, the reduced system becomes O.D.E. if the $G-$
action is of cohomogeneity one. To show the procedure, we consider some
concrete $G-$actions of cohomogeneity one on $CP^n$ and $C^n$ respectively.
By solving the corresponding O.D.E.,we construct infinitely many non-trivial
closed $H-$minimal Lagrangian submanifolds and also non-trivial complete $H-$
minimal Lagrangian submanifolds in $CP^n$ and $C^n$ . Here the $H-$minimal
Lagrangian submanifolds is called nontrivial, if they are not minimal in the
usual sense.

\heading{\bf 2. Symmetry Reduction}
\endheading
\vskip 0.3 true cm
Let $M$ be a connected manifold with a differentiable $G-$action, where $G$
is a compact, connected Lie group. For each $x\in M$ let $G_x$ be the
isotropy subgroup of $x$, and $G(x)\approx G/G_x$ be the orbit of $x$ under $
G$. Two orbits, $G(x)$ and $G(y)$, are said to be of the same type if $G_x$
and $G_y$ are conjugate in $G$. The conjugacy classes of the subgroup $
\{G_x:x\in M\}$ are called the orbit types of the $G-$space $M$. The orbit
types may be partially ordered as follows: 
$$
(H)\succ (K)\Longleftrightarrow \exists g\in G\text{ s.t. }K\supseteq
gHg^{-1}
$$
where $(H)$ denotes the conjugacy class of $H$. We need the following
important result ([MSY]):

\proclaim{Proposition 2.1}
(Principal orbit type) Let $M$ be a connected manifold with a differentiable 
$G-$action. Then there exists a unique orbit type $(H)$ such that $(H)\succ
(K)$ for all orbit types $(K)$ of the action. Moreover, the union of all
orbits of type $(H)$, namely $M^{*}=\{x\in M:G_x\in (H)\}$, is an open,
dense submanifold of $M$.
\endproclaim

Following [MY] we call $\left( H\right) $ in Proposition 2.1 the principal
orbit type of the $G-$space $M$. If $(H^{\prime })\neq (H)$ but $\dim
H^{\prime }=\dim H$, then $(H^{\prime })$ will be called an exceptional
orbit type. All other orbit types will be called singular.

From now on we assume that $M$ is a Kaehler manifold with Kaehler form $
\omega $ and complex structure $J$. Let $G\ $be a compact, connected Lie
group of holomorphic isometries of $M$. Let $g$ be the Lie algebra of $G$,
and $g^{*}$ the dual space of $g$. Then a moment map for the action of $G$ 
on $M$ is a smooth map $\mu:M\rightarrow g^{*}$ such that
\roster
\item "(a)" $(d\mu ,\xi )=i_{\phi (\xi )}\omega $ for all $\xi \in g$ , where $(,)$
denotes the pairing between $g$ and $g^{*}$, and $\phi :g\rightarrow
C^\infty (TM)$ is the infinitesimal action;
\item "(b)" $\mu (kx)=Ad_k^{*}\mu (x)$, $\forall k\in G$ and $x\in M$;
\endroster
where $Ad^{*}$ denotes the coadjoint action.

Let $G$ be an action on $(M,\omega ,J)$ with moment map $\mu $. Let $
Z(g^{*})$ be the centre of $g^{*}$, i.e., the vector subspace of $g^{*}$
fixed by the coadjoint action of $G$. If $c\in Z(g^{*})$, we see from (b)
that $G$ induces an action on the level set $\mu ^{-1}(c)$. Denote by $\pi
:\mu ^{-1}(c)\rightarrow \mu ^{-1}(c)/G$ the natural projection to the
quotient space. Set $\mu ^{*}=$ $\mu |_{M^{*}}$, where $M^{*}$ consists of
principal orbits of the $G-$action. Obviously, 
$$
\mu ^{*}{}^{-1}(c)=\mu ^{-1}(c)\cap M^{*}.
$$
We need the following singular symplectic reduction in Kaehler case(cf.
[SL]) :

\proclaim{Proposition 2.2}
Suppose that $G$ acts on the Kaehler manifold $(M,\omega ,J)$ with moment
map $\mu $ and preserving $J$. If $c\in Z(g^{*})$, then $\mu ^{-1}(c)$ is a
stratified manifold which induces a stratified Kaehler space $\mu
^{-1}(c)/G $ . In particular, $\mu ^{*-1}(c)$ is a manifold and the quotient
space $\mu ^{*}{}^{-1}(c)/G$ inherits a natural Kaehler structure $(
\widetilde{\omega },\widetilde{J})$ such that $\pi ^{*}\widetilde{\omega }
=\omega |_{\mu ^{*-1}(c)}$ and $\widetilde{J}\widetilde{X}=\pi _{*}(JX)$ for
any $\widetilde{X}\in T_p(\mu ^{*}{}^{-1}(c)/G)$ , where $X$ is the
horizontal lift of $\widetilde{X}$.
\endproclaim

\remark{Remark 2.1}
(i) Actually, $\mu ^{*-1}(c)$ is the stratum of $\mu ^{-1}(c)$ corresponding to the principal
orbit type. It is easy to see that $\pi :(\mu ^{*-1}(c),$ $ds^2)\rightarrow
(\mu ^{*-1}(c)/G,\widetilde{ds}^2)$ is a Riemannian submersion, where $ds^2$
is the induced metric from the Kaehler metric of $M$ and $\widetilde{ds}^2$
is the metric determined by $\widetilde{\omega }$;

(ii) If $c\in Z(g^{*})$ is a regular value of $\mu $ and the action of $G$
is free, then $\pi :$ $\mu ^{-1}(c)\rightarrow \mu ^{-1}(c)/G$ is the
well-known Marsden-Weinstein symplectic reduction. In this case, $\mu
^{-1}(c)$ and $\mu ^{-1}(c)/G$ are both smooth manifolds.
\endremark

The following result shows that moment maps are a useful tool for studying
Lagrangian submanifolds with symmetries.

\proclaim{Lemma 2.3}
(cf. [J2]) If $L\subset M$ is a connected $G-$invariant Lagrangian
submanifold, then $L\subset \mu ^{-1}(c)$ for some $c\in Z(g^{*})$.
\endproclaim

\demo{Proof} For $\xi \in g$, we have the vector field $\phi (\xi )$ on $M$. Since 
$\omega |_L\equiv 0$ and $\phi (\xi )$ is tangent to $L$, we have $d\mu
|_{TL}\equiv 0$ by (a). So $\mu $ is constant on $L$. By (b), we see that
the constant $\mu (L)\in Z(g^{*})$. \qed
\enddemo

\proclaim{Lemma 2.4}
Let $i:L\rightarrow M$ be a $G-$invariant Lagrangian submanifold of $M$.
Then $L$ is $H-$minimal if and only if the volume of $L$ is stationary
w.r.t. all compactly supported, $G-$equivariant Hamiltonian variations.
\endproclaim

\demo {Proof} Let $H$ be the mean curvature vector of $L$. Since $H$ depends only
on the immersion $i$, and $i$ is $G-$invariant, we have $k_{*}H=H$ for any $
k\in G$. This implies that the one form $\alpha _H:=H\lrcorner \omega $ and
thus its codifferential $\delta \alpha _H$ are $G-$invariant. Let $\varphi $
be any smooth, $G-$invariant, compactly supported function on $L$. We define
a variation $i_t$, $-\varepsilon <t<\varepsilon $, of the immersion $i$ by : 
$$
i_t(x)=\exp _x(tV)
$$
where $JV=\nabla (\varphi \delta \alpha _H)$, i.e., $\alpha _V=d(\varphi
\delta \alpha _H)$. We choose $\varepsilon >0$ small enough that each $i_t$
is an immersion. Observe that 
$$
\align
k\circ i_t(x) &=k\circ \exp _x(tV) \\
&=\exp _{kx}(tk_{*}V_x) \\
&=\exp _{kx}(tV_{kx}) \\
&=i_t\circ k(x).
\endalign
$$
Hence each $i_t$ is equivariant. By the first variational formula of volume,
we have 
$$
\align
\frac d{dt}|_{t=0}Vol(i_t(L)) &=-\int_L<H,V> \\
&=-\int_L<\alpha _H,\alpha _V> \\
&=-\int_L|\delta \alpha _H|^2\varphi .
\endalign
$$
So, by assumption $\frac d{dt}|_{t=0}Vol(i_t(L))=0$ and the fact that $\varphi$
is arbitrary, we see that $\delta \alpha _H=0$.\qed
\enddemo

We have the generalized Noether Theorem:

\proclaim {Proposition 2.5}
(cf. [Si]) If $F$ is a $G-$invariant Hamiltonian, then the moment map $\mu $
is a conserved quantity for the Hamiltonian flow of $F$.
\endproclaim

Lemma 2.3 shows that any connected $G-$invariant Lagrangian submanifold $L$
is contained in a level set $\mu ^{-1}(c)$ for some $c\in Z(g^{*})$.
Proposition 2.5 implies that the deformation $L_t$ of $L$ by a $G-$invariant
Hamiltonian flow is still contained in the same level set $\mu ^{-1}(c)$.
So we may restrict the equivariant Hamiltonian variational problem in a
fixed level set.

For $c\in Im(\mu ^{*})\cap Z(g^{*})$ , we have a Riemannian submersion $\pi
:(\mu ^{*-1}(c),ds^2)\rightarrow (\mu ^{*-1}(c)/G,\widetilde{ds}^2)$ , whose
fibers are the principal orbits of the $G-$action. Obviously, $\mu ^{*-1}(c)$
and $\mu ^{*-1}(c)/G$ are open dense submanifolds of the stratified spaces $
\mu ^{-1}(c)$ and $\mu ^{-1}(c)/G$ respectively. According to [HsLa], we
define the volume function of the orbits as follows: 
$$
\aligned
V:&\mu ^{*-1}(c)/G \longrightarrow R^{+} \\ 
x &\longmapsto Vol(\pi ^{-1}(x)).
\endaligned \tag{2}
$$
Let $i:L\rightarrow \mu ^{-1}(c)\subset M$ be a $G-$invariant Lagrangian
submanifold, and for simplicity assume that $i(L)\cap \mu ^{*-1}(c)\neq
\emptyset $. (There is no loss of generality in this assumption since $\mu
^{-1}(c)$ may always be replaced by certain natural substrata for which the
assumption holds and to which all subsequent arguments apply.). The
cohomogeneity of $L$ is defined as the integer $\dim L-v$, where $v$ is the
common dimension of the principal orbits. Obviously, if $L$ is of
cohomogeneity $k$, then it project to a map $\widetilde{i}:$ $L/G\rightarrow
\mu ^{-1}(c)/G$ such that $(\widetilde{i}|L^{*}/G):L^{*}/G\rightarrow \mu
^{*-1}(c)/G$ is a $k-$dimensional Lagrangian submanifold of $(\mu ^{*-1}(c)/G,
\widetilde{\omega })$. We will denote $L^{*}/G$ by $\widetilde{L}$. The
Hsiang-Lawson metric on $\mu ^{-1}(c)/G$ is defined as follows (cf. [HsLa]): 
$$
\widetilde{g}_{HL}=V^{2/k}\widetilde{g},  \tag{3}
$$
which goes continuously to zero at the singular boundary.

\proclaim {Theorem 2.6}
Let $i:L\rightarrow M$ be a $G-$invariant Lagrangian submanifold of $M$ with 
$L\subset \mu ^{-1}(c)$ for some $c\in Z(g^{*})$. Then $L$ is a $H-$minimal
Lagrangian submanifold of $M$ if and only if $\widetilde{i}:\widetilde{L}
\rightarrow \mu ^{*-1}(c)/G$ is a $H-$minimal Lagrangian submanifold of $
(\mu ^{*-1}(c)/G,\widetilde{\omega },\widetilde{g}_{HL})$. Furthermore, if $
L $ is of cohomogeneity $k$, then $\widetilde{L}$ is $H-$minimal if and only
if 
$$
\delta _{HL}(V^{4/k}\widehat{H}\lrcorner \widetilde{\omega })=0,  \tag{4}
$$
where $\delta _{HL}$ is the codifferential operator w.r.t. the metric $
\widetilde{g}_{HL}$ and $\widehat{H}$ is the mean curvature vector field of
the Lagrangian submanifold $\widetilde{L}\hookrightarrow (\mu ^{*-1}(c)/G,
\widetilde{\omega },\widetilde{g}_{HL})$.
\endproclaim
\demo {Proof} Denote by $H_G^\infty (T\mu ^{*-1}(c))$ the set of $G-$invariant
horizontal vector fields on $\mu ^{*-1}(c)$ and $C^\infty (T(\mu
^{*-1}(c)/G))$ the set of vector fields on $\mu ^{*-1}(c)/G$. It is easy to
see that 
$$
\align 
\pi _{*}:H_G^\infty (T\mu ^{*-1}(c))&\rightarrow C^\infty (T(\mu ^{*-1}(c)/G))
\\ 
W &\mapsto \widetilde{W}
\endalign
$$
is an bijective correspondence. For any $W,X\in H_G^\infty (T\mu ^{*-1}(c))$
, we have by Proposition 2.2 that 
$$
\align
(W\lrcorner \omega )(X) &=\omega (W,X) \\
&=(\pi ^{*}\widetilde{\omega })(W,X) \\
&=\widetilde{\omega }(\widetilde{W},\widetilde{X}) \\
&=(\widetilde{W}\lrcorner \widetilde{\omega })(\widetilde{X})
\endalign
$$
where $\widetilde{W}=\pi _{*}W$ and $\widetilde{X}=\pi _{*}X$. Obviously, $
W|_L$ is a Hamiltonian vector field along $L$ w.r.t. $\omega $ if and only
if $\widetilde{W}|_{\widetilde{L}}$ is a Hamiltonian field along $\widetilde{
L}$ w.r.t. $\widetilde{\omega }$. On the other hand, (3) implies that: 
$$
Vol(L,g)=Vol(\widetilde{L},\widetilde{g}_{HL}).
$$
Thus the first part of the Theorem follows immediately from Lemma 2.4.

Let $f$ be an arbitrary compactly supported function on $\widetilde{L}$ ,
which determines a Hamiltonian normal vector field $\widetilde{W}$ along $
\widetilde{L}$, i.e., $\widetilde{W}\lrcorner \widetilde{\omega }=df$. By
the first variational formula we have

$$
\align
\frac{dVol((\widetilde{L}_t,\widetilde{g}_{HL}))}{dt}|_{t=0}&=-\int_{
\widetilde{L}}<\widehat{H},\widetilde{W}>_{\widetilde{g}_{HL}}dVol_{HL} \\ 
&=-\int_{\widetilde{L}}<J\widehat{H},J\widetilde{W}>_{\widetilde{g}_{HL}}dVol_{HL} \\ 
&=-\int_{\widetilde{L}}V^{4/k}<\widehat{H}\lrcorner \widetilde{\omega },
\widetilde{W}\lrcorner \widetilde{\omega }>_{\widetilde{g}_{HL}}dVol_{HL} \tag {5}\\ 
&=-\int_{\widetilde{L}}V^{4/k}<\widehat{H}\lrcorner \widetilde{\omega },df>_{
\widetilde{g}_{HL}}dVol_{HL} \\ 
&=-\int_{\widetilde{L}}\delta _{HL}(V^{4/k}\widehat{H}\lrcorner \widetilde{
\omega })fdVol_{HL}.
\endalign 
$$
Then (4) follows immediately from (5). \qed
\enddemo

\proclaim {Corollary 2.7}
Under the assumption of Theorem2.6, if $\widetilde{L}$ is minimal in $
(\mu ^{-1}(c)/G,\widetilde{g}_{HL})$ then $L$ is $H-$minimal in $M$.
\endproclaim

Theorem 2.6 is interesting and particular simple when all the $G-$orbits are
isometric. In this case, the volume function of orbits is constant, and thus
the metrics $\widetilde{g},\widetilde{g}_{HL}$ are equivalent.

\proclaim {Corollary 2.8}
If all $G-$orbits in $\mu ^{-1}(c)\,\,$are mutually isometric , then $
L\subset \mu ^{-1}(c)$ is a $G-$invariant $H-$minimal submanifold in $M$ if
and only if $\pi (L)$ is a $H-$minimal submanifold in $(\mu ^{-1}(c)/G,
\widetilde{\omega },\widetilde{g})$.
\endproclaim

Let us consider an important special case of Corollary 2.8. Define a $S^1-$
action on $C^n$ by
$$
e^{i\theta }(z_1,...,z_n)=(e^{i\theta }z_1,...,e^{i\theta }z_n).  \tag{6}
$$
This is a Hamiltonian action with moment map 
$$
\mu (z)=-\frac i2|z|^2.  \tag{7}
$$
Its level set at a value $-\frac i2t$ is $S^{2n-1}(\sqrt{t})$. The
symplectic reduction at a regular value $-\frac i2t$ ( $t>0$) gives a
fibration $\mu ^{-1}($ $-\frac i2t)\rightarrow \mu ^{-1}(-\frac i2t)/S^1$.
In particular, we have the well-known Hopf-fibration $\pi
:S^{2n-1}\rightarrow CP^{n-1}$ by taking $t=1$. Since all $S^1-$orbits in $
\mu ^{-1}($ $-\frac i2)=S^{2n-1}$ are isometric, we obtain immediately from
Corollary 2.8 the following result:

\proclaim {Theorem 2.9}
Let $\pi :S^{2n-1}\rightarrow CP^{n-1}$ be the Hopf fibration. Let $
\widetilde{L}^{n-1}\hookrightarrow CP^{n-1}$ be a Lagrangian submanifold and $
L^n=\pi ^{-1}(\widetilde{L}^{n-1})$ the inverse image of $\widetilde{L}
^{n-1} $ by the Hopf projection. Then $L^n$ is a $H-$minimal Lagrangian
submanifold in $C^n$ if and only if $\widetilde{L}^{n-1}$ is a $H-$minimal
Lagrangian submanifold in $CP^{n-1}$.
\endproclaim

\remark {Remark 2.2}
(i) In [Oh2], it was proved that the inverse image $\pi ^{-1}(\widetilde{L}
^{n-1})$ is a $H-$minimal Lagrangian submanifold in $C^n$ provided that $
\widetilde{L}^{n-1}$ is a usual minimal Lagrangian submanifold in $CP^{n-1}$
. So Theorem 2.9 generalizes Oh's result; (ii) In [HR2], the authors
constructed $H-$minimal Lagrangian tori in $CP^2$, which are not minimal, by
integrable system method. By applying Theorem 2.9 to their examples, we can
get a large number of non-trivial $H-$minimal $T^3$ in $C^3$; (iii)
Corollary 2.7 may also be regarded as a generalization of Oh's result in
another direction.
\endremark

If the Lagrangian submanifold is of cohomogeneity one, we may simplify the
equation (4) as follows:

\proclaim {Corollary 2.10}
Let $L\hookrightarrow \mu ^{-1}(c)\subset M$ be a $G-$invariant Lagrangian
submanifold of cohomogeneity one, then $L$ is $H-$minimal if and only if 
$$
V^2k_{\widetilde{L}}=K\text{.}  \tag{8}
$$
where $k_{\widetilde{L}}$ is the mean curvature of the curve $\widetilde{L}$ in $(\mu ^{*-1}(c)/G,\widetilde{g}
_{HL})$ and where $K$ is any constant.
\endproclaim

\demo {Proof} Let $\widehat{e}_1$ be the unit tangent vector field of the curve $
\widetilde{L}$ with respec to the Hsiang-Lawson metric. From (4), we have 
$$
\align
0 &=\widehat{e}_1(V^{2}\widetilde{g}_{HL}(J\widehat{H},\widehat{e}_1)) \\
&=-\widehat{e}_1(V^{2}\widetilde{g}_{HL}(\widehat{H},J\widehat{e}_1)) \\
&=\widehat{e}_1(V^{2}k_{\widetilde{L}}),
\endalign
$$
i.e., $V^2k_{\widetilde{L}}=$const.\qed
\enddemo

\remark {Remark 2.3}
Corollary 2.10 reduces the third order P.D.E (1) to a second order O.D.E.
(8) with a constant $K$.
\endremark

In the remaining two sections, we will use Corollary 2.10 to construct
Hamiltonian minimal submanifolds of cohomogeneity one in $CP^n$ and $C^n$
respectively.

\heading {\bf 3. Hamiltonian minimal Lagrangian submanifolds in $CP^n$}
\endheading
\vskip 0.3 true cm
\subhead {3.1. $SO_n-$invariant $H-$minimal Lagrangian submanifolds}
\endsubhead
\vskip 0.3 true cm
Let $G=$ $SO(n)$, that can be regarded as a subgroup of $PU(n+1)=U(n+1)/S^1$
in the natural way. The group $G$ acts on $(CP^n,\omega _{FS})$ by 
$$
A\cdot [z]=[z_0:\widehat{z}_1:\cdots :\widehat{z}_n]  \tag{9}
$$
where $z=(z_0,z_1,...,z_n)$ and $(\widehat{z}_1,....,\widehat{z}
_n)^t=A(z_1,...,z_n)^t$. This is a Hamiltonian action on $CP^n$, whose
moment map is given by 
$$
\mu ([z])=\frac 1{|z|^2}(Im(z_1\overline{z}_2),...,Im(z_1\overline{z}
_n),Im(z_2\overline{z}_3),...,Im(z_2\overline{z}_n),...,
Im(z_{n-1}\overline{z}_n)).
\tag{10}
$$
As $Z(g^{*})=\{0\}$, any $G-$invariant connected Lagrangian submanifold lies
in $\mu ^{-1}(0)$. All points in $\mu ^{-1}(0)$ may be written as 
$$
\lbrack x_0:\lambda x_1:\lambda x_2:\cdots :\lambda x_n]
$$
where $\lambda \in C$ and $x_0,x_1,...,x_n$ are real, and normalized so that 
$\sum_{\alpha =1}^nx_\alpha ^2=1$ and $x_0^2+|\lambda |^2=1$. Therefore the
orbits of $G$ in $\mu ^{-1}(0)$ are $O_\lambda $ for $\lambda \in C$, where 
$$
O_\lambda =\{[x_0:\lambda x_1:\lambda x_2:\cdots :\lambda x_n]:x_\alpha \in
R,x_0^2+|\lambda |^2=1,\sum_{i=1}^nx_i^2=1\}.
$$
The orbit space $\mu ^{-1}(0)/G$ may be parameterized as 
$$
\mu ^{-1}(0)/G=\{[\sqrt{1-r^2},re^{i\theta },0,...,0]:0\leq r\leq 1,\theta
\in [0,2\pi )\}.
$$
Note that $r=0$ and $r=1$ correspond to a singular orbit and an exceptional
orbit respectively. Each orbit in $\mu ^{*-1}(0)/G$ has the following unique
representative element in $S^{2n+1}(1)$ 
$$
\align
F(r,\theta )=(\sqrt{1-r^2},0, r\cos \theta ,&r\sin \theta ,0,...,0)\in
S^{2n+1}(1), \tag{11} \\ 
0<r<1,\quad &0\leq \theta <2\pi .
\endalign 
$$
At $F(r,\theta )$, the unit vertical vector $\eta $ of the Hopf fibration $
\pi _H:S^{2n+1}(1)\rightarrow CP^n$ is given by 
$$
\eta (r,\theta )=(0,\sqrt{1-r^2},-r\sin \theta ,r\cos \theta ,0,...,0). 
\tag{12}
$$
To determine the tangent space of the $G-$orbit $O_\lambda $, we consider
any tangent vector $X=(0,v_2...,v_n)\in $ $T_pS^{n-1}(1)$ at $
p=(1,0,...,0)\in S^{n-1}(1)$. Set 
$$
\xi _X=(0,0,0,0,rv_2\cos \theta ,rv_2\sin \theta ,...,rv_n\cos \theta
,rv_n\sin \theta )\in R^{2n+2}=C^{n+1}.  \tag{13}
$$
Obviously $<\xi _X,\eta >=0$, and thus $Span\{(\pi _H)_{*}\xi _X$ $:X\in
T_pS^{n-1}(1)\}$ is just the tangent space of the $G-$orbit at the
corresponding point. From (11), we have 
$$
dF(\frac \partial {\partial r})=(\frac{-r}{\sqrt{1-r^2}},0,\cos \theta ,\sin
\theta ,0,...,0)
$$
and 
$$
dF(\frac \partial {\partial \theta })=(0,0,-r\sin \theta ,r\cos \theta
,0,....,0).
$$
Obviously $<dF(\frac \partial {\partial r}),\eta >=<dF(\frac \partial
{\partial r}),\xi _X>=0$ . So, 
$$
|\pi _{*}\pi _{H*}(dF(\frac \partial {\partial r})|^2=\frac 1{1-r^2}.
$$
The horizontal component (w.r.t. $\pi _H$) of $dF(\frac \partial {\partial
\theta })$ is given by: 
$$
F_\theta :=dF(\frac \partial {\partial \theta })-r^2\eta \text{.}  \tag{14}
$$
From (12) , (13) and (14), we see that $<F_\theta ,\xi _X>=0$ . Then 
$$
|\pi _{*}\pi _{H*}(F_\theta )|^2=r^2-r^4. 
$$
Also $<dF(\frac \partial {\partial r}),F_\theta >=0$. Hence the induced
metric on the orbit space $\mu ^{-1}(0)/G$ is given by 
$$
\widetilde{g}=\frac 1{1-r^2}dr^2+r^2(1-r^2)d\theta ^2. \tag{15}
$$
Up to a constant the volume function of the orbits is $r^{n-1}$. Therefore
the Hsiang-Lawson metric on $\mu ^{-1}(0)/G$ is given by 
$$
\widetilde{g}_{HL}=r^{2n-2}[\frac 1{1-r^2}dr^2+r^2(1-r^2)d\theta ^2]
$$
If we set $r=\sin \varphi $, $\widetilde{g}_{HL}$ can be expressed as 
$$
\aligned 
\widetilde{g}_{HL}=&\sin ^{2n-2}\varphi [d\varphi ^2+\sin ^2\varphi \cos
^2\varphi d\theta ^2],\\ 
&0\leq \varphi \leq \pi /2,\ 0\leq \theta <2\pi .
\endaligned \tag{16}
$$
Here $\theta $ is the rotational parameter, and $\varphi $ is the radial
parameter. Note that $\varphi =0$ and $\varphi =\frac \pi 2$ correspond to
singular points on $\mu ^{-1}(0)/G$.

Observe that the metric is invariant under the rotation in $\theta $ and the
reflection $\theta _0+\theta \rightarrow \theta _0-\theta $ for any $\theta
_0$. So $\theta \equiv const.$ are all geodesics on $(\mu ^{-1}(0)/G)$,
whose inverse images are mutually congruent in $CP^n$. The congruence class
corresponds to the totally geodesic Lagrangian immersion $S^n\rightarrow
RP^n\subset CP^n$. We are not interested in this case.

Let $\varphi (\theta )$ be any curve in $\mu ^{*-1}(0)/G$, where $\theta $
is now allowed to vary over all real numbers. The unit tangent vector field
and the normal vector field of $(\theta ,\varphi (\theta ))$ are given
respectively by 
$$
e=\frac 1{\sin ^{n-1}\varphi \sqrt{(\varphi ^{\prime })^2+\sin ^2\varphi
\cos ^2\varphi }}(\varphi ^{\prime }\frac \partial {\partial \varphi }+\frac
\partial {\partial \theta })
$$
and 
$$
n=\frac 1{\sin ^n\varphi \cos \varphi \sqrt{(\varphi ^{\prime })^2+\sin
^2\varphi \cos ^2\varphi }}(-\sin ^2\varphi \cos ^2\varphi \frac \partial
{\partial \varphi }+\varphi ^{\prime }\frac \partial {\partial \theta }).
$$
For any variation $\varphi +s\eta $ of $\varphi $, we have the corresponding
variation vector field 
$$
\xi =\eta \frac \partial {\partial \varphi }.  \tag{17}
$$
So we get 
$$
<\xi ,n>=-\frac{\eta \sin ^n\varphi \cos \varphi }{\sqrt{(\varphi ^{\prime
})^2+\sin ^2\varphi \cos ^2\varphi }}.  \tag{18}
$$
By definition, $k_{\widetilde{L}}=<\nabla _e^{HL}e,n>_{HL}$, where $\nabla
^{HL}$ denotes the Levi-Civita connection of the metric $\widetilde{g}_{HL}$
. From Corollary 2.10, we know that the $H-$minimal equation for $\widetilde{
L}$ is 
$$
k_{\widetilde{L}}\sin ^{2n-2}\varphi =K.  \tag{19}
$$
It is easy to see from the first variation formula of the arc length w.r.t.
the variation (17) that (19) is the Euler-Lagrange equation for the
following functional: 
$$
J=\int [\sin ^{n-1}\varphi \sqrt{\sin ^2\varphi \cos ^2\varphi +(\varphi
^{\prime })^2}-\frac K2\sin ^2\varphi ]d\theta .  \tag{20}
$$
By a direct computation from (20), we get the E-L equation for $\varphi
(\theta )$: 
$$
\aligned
\frac{\sin ^{n-1}\varphi }{(\sqrt{\sin ^2\varphi \cos ^2\varphi +\varphi
^{\prime 2}})^3}&\{-\varphi ^{\prime \prime }\sin \varphi \cos \varphi
+[(n+1)\cos ^2\varphi -2\sin ^2\varphi ](\varphi ^{\prime })^2 \\ 
+&\sin ^2\varphi \cos ^2\varphi [n\cos ^2\varphi -\sin ^2\varphi ]\}=K,
\endaligned \tag{21}
$$
where $K$ is a constant. We will assume $K\neq 0$, because this condition
keeps the $H-$minimal submanifolds from being minimal submanifolds. Set 
$$
L(\theta ,\varphi ,\varphi ^{\prime })=\sin ^{n-1}\varphi \sqrt{\sin
^2\varphi \cos ^2\varphi +(\varphi ^{\prime })^2}-\frac K2\sin ^2\varphi .
$$
We perform a Legendre transformation 
$$
p=L_{\varphi ^{\prime }}=\frac{\varphi ^{\prime }\sin ^{n-1}\varphi }{\sqrt{%
\sin ^2\varphi \cos ^2\varphi +(\varphi ^{\prime })^2}}.
$$
The Hamiltonian $H$ of the equation (21) is defined as 
$$
\aligned
H(\theta ,\varphi ,p)&=\varphi ^{\prime }p-L \\ 
&=-\sin \varphi \cos \varphi \sqrt{\sin ^{2n-2}\varphi -p^2}+\frac K2\sin
^2\varphi .
\endaligned  \tag{22}
$$
Note that $H$ does not depend explicitly on the variable $\theta $. So $H$
is a constant of motion, i.e., constant along any solution of the equation
(cf. [JL]). It follows from (22) that 
$$
\frac{\sin ^{n+1}\varphi \cos ^2\varphi }{\sqrt{\sin ^2\varphi \cos
^2\varphi +(\varphi ^{\prime })^2}}=\lambda +\frac K2\sin ^2\varphi , 
\tag{23}
$$
where $\lambda $ and $K$ are constants.

We will solve the ODE (21) by considering the following initial conditions 
$$
\aligned
\varphi (0)&=a\in (0,\frac \pi 2), \\ 
\varphi ^{\prime }(0)&=b.
\endaligned \tag{24}
$$
So $\lambda $ is determined by the initial conditions from (23) as follows 
$$
\lambda =\frac{\sin ^{n+1}a\cos ^2a}{\sqrt{\sin ^2a\cos ^2a+b^2}}-\frac
K2\sin ^2a.  \tag{25}
$$
At least the ODE (23) with the initial values (24) can be solved locally.
Any such a solution gives a (local) H-minimal Lagrangian submanifold in $CP^n
$. We now give the initial values $(a,b)$ to ensure the global existence of
the solution.

\proclaim {Lemma 3.1}
For any initial values $(a,b)$ with $\lambda \notin [0,-\frac K2]$ or $
[-\frac K2,0]$ according to $K<0$ or $K>0$, there is a unique global
solution of (21) satisfying the initial conditions (24).
\endproclaim

\demo {Proof} Set $\psi =\varphi ^{\prime }$. Then we may rewrite the ODE (21) as
an ODE system of first order for $(\theta ,\varphi ,\psi )$ on the domain $
(-\infty ,\infty )\times (0,\frac \pi 2)\times (-\infty ,\infty )$. Under
the hypothesis in the Lemma, we see from (23) that there exists no finite
value $\theta _0\in R$ such that $\varphi (\theta )\rightarrow 0$ , $\frac
\pi 2$ or $\psi (\theta )\rightarrow \infty $ as $\theta \rightarrow \theta
_0$. Hence the local solution may be extended to a global solution.\qed
\enddemo

Now we hope to determine the initial values for which the corresponding $
\lambda $ satisfies the condition of Lemma 3.1. First, if $K<0$, it is easy
to see from (25) that $\lambda \notin [0,-\frac K2]$ is equivalent to 
$$
\frac{\sin ^{n+1}a}{\sqrt{\sin ^2a\cos ^2a+b^2}}>-\frac K2.  \tag{26}
$$
Similarly, if $K>0$, the condition $\lambda \notin [-\frac K2,0$ $]$ is
equivalent to 
$$
\frac{\sin ^{n-1}a\cos ^2a}{\sqrt{\sin ^2a\cos ^2a+b^2}}>\frac K2.  \tag{27}
$$
Obviously we can always find initial values $(a,b)$ such that (26) or (27)
is satisfied, provided that $|K|$ is small enough. Hence the ODE (21) has a
global solution for any initial values (24) which satisfy (26) or (27)
according to $K<0$ or $K>0$. From the proof of Lemma 3.1, we see that $
\inf_{(-\infty ,\infty )}\sin ^2\varphi (\theta )\cos ^2\varphi (\theta )=B>0
$ for the global solution. The length of $\varphi (\theta )$ with respect to 
$\widetilde{g}$ given by (15) is 
$$
\align
L(\varphi ) &=\int_{-\infty }^\infty \sqrt{(\varphi ^{\prime })^2+\sin
^2\varphi \cos ^2\varphi }d\theta \\
&\geq \sqrt{B\;}\int_{-\infty }^\infty d\theta \\
&=\infty ,\endalign
$$
i.e., the solution curve has infinite length. As a Riemannian manifold of
one dimension, the solution curve is complete. Since the fibres of the
projection $\pi :\mu ^{-1}(0)\rightarrow \mu ^{-1}(0)/G$ are compact, it is
easy to prove that the corresponding $H-$minimal submanifold $\pi ^{-1}(
\widetilde{L})$ is complete as a metric space. So it is complete via
Hopf-Rinow Theorem. We have proved

\proclaim {Theorem 3.2}
There exist infinitely many non-trivial complete $H-$minimal Lagrangian
immersions of $R^1\times S^{n-1}$ into $CP^n$.
\endproclaim

If $\varphi =\varphi (\theta )$ corresponds to a closed curve, then there is
some point, which we may assume is $\theta =0$, at which $\varphi $ assumes
a maximum or minimum. Hence we consider the following initial conditions 
$$
\aligned
\varphi (0)&=a \\ 
\varphi ^{\prime }(0)&=0
\endaligned  \tag{28}
$$
for $a\in (0,\frac \pi 2)$.

According to Lemma 3.1, we will choose the initial value $a$ such that 
$$
\align
a\in I_K:&=\{x\in (0,\frac \pi 2)|\sin ^nx\sec x>-\frac K2\text{ for }K<0 \\ 
&\text{or }\sin ^{n-2}x\cos x>\frac K2\text{ for }K>0\}.
\endalign
$$
So we get a global solution $\varphi $ for such initial values. If $\varphi
^{\prime \prime }(0)=0$, then $\varphi \equiv const.$ and the constant
solution can be determined from (21).

Without lose of generality, we assume that $\varphi ^{\prime \prime }(0)<0$.
From (23), we have 
$$
\frac{d\varphi }{d\theta }=\pm \frac{\sin \varphi \cos \varphi \sqrt{\sin
^{2n}\varphi \cos ^2\varphi -(\lambda +\frac K2\sin ^2\varphi )^2}}{\lambda
+\frac K2\sin ^2\varphi }  \tag{29}
$$
and thus 
$$
\theta =\pm \int_a^{\varphi (\theta )}\frac{(\lambda +\frac K2\sin ^2\varphi
)d\varphi }{\sin \varphi \cos \varphi \sqrt{\sin ^{2n}\varphi \cos ^2\varphi
-(\lambda +\frac K2\sin ^2\varphi )^2}}.  \tag{30}
$$
Set $f(x):=\sin ^{2n}x\cos ^2x-(\lambda +\frac K2\sin ^2x)^2$. We get from
(28) and (29) that $f(a)=0$ and 
$$
\frac{d^2\varphi }{d\theta ^2}(0)=\frac{f^{\prime }(a)}{2\sin ^{2n-2}a}<0
$$
i.e., $f^{\prime }(a)<0$. So  $f(a-\varepsilon )>0$ for small $\varepsilon >0
$. On the other hand, $f(0)=-\lambda ^2<0$. Thus there is $b\in (0,a)$
such that $f(b)=0$. Set $\widehat{b}=\max \{b:f(b)=0,0<b<a\}$ and 
$$
\Omega _a=\min \{\theta |\varphi (\theta )=\widehat{b},\theta \in (0,+\infty
)\}.  \tag{31}
$$
Then $\varphi $ is a decreasing function on $[0,\Omega _a]$ and $\varphi
(\Omega _a)=\widehat{b}$. By (30), $\Omega _a$ is given by 
$$
\Omega _a=-\int_a^{\widehat{b}}\frac{(\lambda +\frac K2\sin ^2\varphi
)d\varphi }{\sin \varphi \cos \varphi \sqrt{\sin ^{2n}\varphi \cos ^2\varphi
-(\lambda +\frac K2\sin ^2\varphi )^2}}.  \tag{32}
$$
We also note that the solution of (21) is invariant under the reflection 
$$
\theta _0+\theta \rightarrow \theta _0-\theta 
$$
for any $\theta _0$. By reflection at points $\{0,\pm n\Omega _a;n=1,2,...\}$
, we get a global solution $\widetilde{\varphi }(\theta )$ on $(-\infty
,\infty )$ with period $\Omega _a$. By uniqueness Theorem of ODE, $\varphi
\equiv \widetilde{\varphi }$ . Obviously, the solution curve $\widetilde{L}$
is closed if and only if $\Omega _a$ is a rational multiple of $\pi $. Since 
$\Omega _a$ is a non-constant continuous function of $a$, we may obtain
countable many such closed curves. Set 
$$
\left. A_n(K)=\{a\in I_K\ |\ \Omega _a/\pi \text{ is rational}\}.\right. 
$$
Let $L_a$ denote the inverse image in $CP^n$ of the closed solution curve $
\varphi (\theta )$ with the initial value $a\in A_n(K)$. Then we have

\proclaim {Theorem 3.3}
There exist countable infinite non-trivial closed $H-$minimal Lagrangian
submanifolds $\{L_a\}_{a\in A_n(K)}$ in $CP^n$ , which are invariant under
the $SO(n)-$ action.
\endproclaim

\remark {Remark 3.1}
These submanifolds $\{L_a\}$ are immersions of $S^1\times S^{n-1}$ in $CP^n$.
\endremark

\vskip 0.4 true cm
\subhead {3.2. $T^{n-1}-$invariant $H-$minimal Lagrangian submanifolds}
\endsubhead
\vskip 0.3 true cm
We only consider the following simple $T^{n-1}-$
action on $CP^n$: 
$$
(e^{i\theta _1},...,e^{i\theta _{n-1}})\cdot [z]=[z_0:e^{i\theta
_1}z_1:\cdots :e^{i\theta _{n-1}}z_{n-1}:z_n]  \tag{33}
$$
whose moment map is 
$$
\mu ([z])=-\frac i{2|z|^2}(|z_1|^2,|z_2|^2,....,|z_{n-1}|^2).  \tag{34}
$$
As $G=T^{n-1}$ is Abelian, $Z(g^{*})=g^{*}$. Choose $c_i\in R$ such that $
c_j>0$ $($ $j=1,...,n-1)$ and $\sum_{j=1}^{n-1}c_j<1$. Set $c=-\frac
i2(c_1,c_2,...,c_{n-1})$. Then we have the level set 
$$
\mu ^{-1}(c)=\{[z]\in CP^n:z\in S^{2n+1},|z_j|^2=c_j,1\leq j\leq n-1\}.
$$
We may parametrize the orbit space $\mu ^{-1}(c)/G$ as follows: 
$$
(r,\theta)\overset{F}\to {\longrightarrow }(\sqrt{1-\sum_{j=1}^{n-1}c_j-r^2},0,
\sqrt{c_1},0,\sqrt{c_2},0,...,\sqrt{c_{n-1}},0,r\cos \theta ,r\sin \theta).\tag{35}
$$
From (35), we have 
$$
\aligned
dF(\frac \partial {\partial r})&=(-\frac r{\sqrt{1-\sum_{j=1}^{n-1}c_j-r^2}},
\underset{2(n-1)}\to {\underbrace{0,0,....,0}},\cos \theta ,\sin \theta ), \\ 
dF(\frac \partial {\partial \theta })&=(\underset{2n}\to {\underbrace{0,.....,0}}
,-r\sin \theta ,r\cos \theta ).\endaligned \tag{36}
$$
Set 
$$
\xi _j=(\underset{2j+1}\to {\underbrace{0,...,0},}\sqrt{c_j},0,...,0),1\leq
j\leq n-1.  \tag{37}
$$
Then $d\pi _H(\xi _j)$ ($j=1,...,n-1$) span the tangent space of the $G-$
orbit at the corresponding point, where $\pi _H:S^{2n+1}\rightarrow CP^n$ is
the Hopf fibration. It is easy to see from (36), (37) that the volume
function of the orbits is constant. Using the similar method as in section
3.1, we can also get from (36), (37) the Hsiang-Lawson metric on $\mu
^{-1}(c)/G$ (up to a constant): 
$$
\widetilde{g}_{HL}=\frac 1{\delta -r^2}dr^2+\frac{r^2(\delta -r^2)}{\delta ^2
}d\theta ^2, \tag{38}
$$
where $\delta =1-\sum_{j=1}^{n-1}c_j$. If we introduce $r=\sqrt{\delta }\sin
\varphi $, then 
$$
\widetilde{g}_{HL}=d\varphi ^2+\sin ^2\varphi \cos ^2\varphi d\theta ^2,\tag{39}
$$
where $0\leq \varphi <\frac \pi 2$ and $0\leq \theta <2\pi $. From Corollary
2.10, we know that $H-$minimal equation for $\widetilde{L}$ is 
$$
k_{\widetilde{L}}=K,  \tag{40}
$$
where $K$ is a constant.

Write $\widetilde{L}$ as $\varphi =\varphi (\theta )$. Similar to the
discussion in section 3.1, we see that (40) is just the Euler-Lagrange
equation of the following functional 
$$
J=\int_{\theta _1}^{\theta _2}(\sqrt{(\varphi ^{\prime })^2+\sin ^2\varphi
\cos ^2\varphi }-\frac K2\sin ^2\varphi )d\theta   \tag{41}
$$
Since (41) is only a special case of (20), we obtain the following:

\proclaim {Theorem 3.4}\newline
(i) There exist infinitely many non-trivial complete $H-$minimal Lagrangian
immersions of $R^1\times T^{n-1}$ into $CP^n$;\newline
(ii) There exist countable infinite non-trivial $H-$minimal
Lagrangian immersions of $T^n\ $into $CP^n$.
\endproclaim

\heading {\bf 4. Hamiltonian minimal Lagrangian submanifolds in $
C^{n+1} $}
\endheading
\vskip 0.3 true cm
\subhead {4.1. Inverse images of the Hopf map}
\endsubhead
\vskip 0.3 true cm
From Theorem 2.9, Theorem 3.3 and Theorem 3.4, we immediately have the
following:

\proclaim {Theorem 4.1}
There exist infinitely many non-trivial $H-$minimal Lagrangian immersions of 
$R^1\times S^1\times S^{n-1},$ $S^1\times S^1\times S^{n-1},$ $R^1\times T^n$
and $T^{n+1}$ into $C^{n+1}$.
\endproclaim

In the rest of this paper, we will consider two Hamiltonian actions on $%
C^{n+1}$, which were used to construct special Lagrangian submanifolds in $%
C^{n+1}$ by Havey and Lawson [HL](cf. also [J2]). We now use them to
construct some new non-trivial complete $H-$minimal Lagrangian submanifolds
in $C^{n+1}$.
\vskip 0.3 true cm
\subhead {4.2. $SO_{n+1}-$invariant $H-$minimal Lagrangian submanifolds}
\endsubhead
\vskip 0.3 true cm
Let $G=SO(n+1)\subset SU(n+1)$, which acts on $C^{n+1}$ $(n>1)$ in the
following way: 
$$
\gamma \cdot z=(\gamma x,\gamma y),z\in C^{n+1},\gamma \in SO(n+1)  \tag{42}
$$
where we write $z=x+iy$. Then the moment map of the action is 
$$
\mu (z_1,...,z_{n+1})=(Im(z_1\overline{z}_2),...,Im(z_1\overline{z}
_{n+1}),Im(z_2\overline{z}_3),...,Im(z_2\overline{z}_{n+1}),...,Im(z_n\overline{z}_{n+1}))
$$
As $Z(g^{*})=\{0\}$, any $G-$invariant connected Lagrangian submanifold is
contained in $\mu ^{-1}(0)$. Obviously any point of $\mu ^{-1}(0)$ may be written as $(\lambda
x_1,...,\lambda x_{n+1})$ with $\lambda \in C$ , $x_1,...,x_{n+1} \in R$ and 
$x_1^2+\cdots +x_{n+1}^2=1$. So a $G-$orbit in $\mu ^{-1}(0)$ is 
$$
O_\lambda =\{(\lambda x_1,\cdots ,\lambda x_{n+1}):x_j\in R,x_1^2+\cdots
+x_{n+1}^2=1\}.
$$
Clearly $O_0$ is a point, and $O_\lambda =O_{-\lambda }\cong S^n$ if $
\lambda \neq 0$. So the orbit space $\mu ^{-1}(0)/G$ is 
$$
\mu ^{-1}(0)/G=\{[(\lambda ,0,...,0)]:\lambda \in C\}.
$$
which has the following parametrization 
$$
(r,\theta )\rightarrow [(r\cos \theta ,r\sin \theta ,0,...,0)].
$$
It is easy to see that the induced metric $\widetilde{g}$ on the orbit space
is flat, i.e., 
$$
\widetilde{g}=dr^2+r^2d\theta ^2.  \tag{43}
$$
The volume of the orbit space at $\lambda $ is (up to a constant): 
$$
V=r^n.
$$
So, the Hsiang-Lawson metric is given by 
$$
\widetilde{g}_{HL}=r^{2n}(dr^2+r^2d\theta ^2).
$$
Obviously, $\theta \equiv const.$ is a geodesic $(\mu ^{-1}(0)/G,\widetilde{g
}_{HL})$, whose inverse image in $C^{n+1}$ is a $(n+1)-$dim Lagrangian plane
passing through the origin of $C^{n+1}$. Now we allow $\theta $ to vary over
all real number and write the curve on the orbit space as $r=r(\theta )$.
Similar to the above discussion, we see that the Hamiltonian-minimal equation for $
r(\theta )$ is 
$$
\frac 1{(r^2+\overset{.}\to {r}^2)^{3/2}}(-r\overset{..}\to {r}+2\overset{.}\to {r}
^2+r^2)+\frac n{\sqrt{r^2+\overset{.}\to {r}^2}}=\frac K{r^n},  \tag{44}
$$
which is a critical point of the following functional 
$$
J=\int (r^n\sqrt{r^2+\overset{.}\to {r}^2}-\frac K2r^2)d\theta ,
$$
where $K$ is a nonzero constant. Since $r=0$ corresponds to a singular point
of the orbit space, we will solve the equation (44) on $(\mu
^{-1}(0)/G)-\{r=0\}$. Set 
$$
L(\theta ,r,\overset{.}\to {r})=r^n\sqrt{r^2+\overset{.}\to {r}^2}-\frac K2r^2 
\tag{45}
$$
and 
$$
p=L_{\overset{.}\to {r}}=\frac{r^n\overset{.}\to {r}}{\sqrt{r^2+\overset{.}\to {r}^2}}. 
\tag{46}
$$
We get the Hamiltonian for the equation: 
$$
H(\theta ,r,p)=-r\sqrt{r^{2n}-p^2}+\frac K2r^2.  \tag{47}
$$
Since $H(\theta ,r,p)$ does not depend on $\theta $ explicitly, it must be a
constant of the motion. It follows from (46) and (47) that 
$$
\frac{r^{n+2}}{\sqrt{r^2+\overset{.}\to {r}^2}}=\lambda +\frac K2r^2,  \tag{48}
$$
where $\lambda ,K$ are constants. For initial values 
$$
\aligned 
r(0)&=a>0, \\ 
r^{\prime }(0)&=b,
\endaligned \tag{49}
$$
we have 
$$
\lambda =\frac{a^{n+2}}{\sqrt{a^2+b^2}}-\frac K2a^2. \tag{50}
$$
From (48), we see that there exists no point $\theta _0$ such that $r(\theta
)\rightarrow 0$ as $\theta \rightarrow \theta _0$, provided that $\lambda
\neq 0$. Any local solution of (44) with the initial values (49) corresponds to
a local $H-$minimal submanifolds in $C^{n+1}$.

In following, we always consider the initial values with $\lambda \neq 0$.
From (48), we have 
$$
\frac{dr}{d\theta }=\pm \frac{r\sqrt{r^{2(n+1)}-(\lambda +\frac K2r^2)^2}}{
\lambda +\frac K2r^2}.  \tag{51}
$$

\proclaim {Lemma 4.2}
If the solution curve $r=r(\theta )$ of (44) is closed, then $r\equiv (\frac
K{n+1})^{1/(n-1)}$ ($K>0$).
\endproclaim

\demo {Proof} Since $r=r(\theta )$ is closed, then there are two points $\theta _1$
and $\theta _2$, at which $r$ assumes the maximum and minimum respectively.
Set $A_i=r(\theta _i),i=1,2$. So, we have 
$$
\aligned
r^{\prime }(\theta _i)& =0,\text{ }i=1,2,\\ 
r^{\prime \prime }(\theta _1)&\leq 0,\ r^{\prime \prime }(\theta _2)\geq 0.
\endaligned \tag 52
$$
Obviously, $\lambda =A_i^{n+1}-\frac K2A_i^2$. From (51), and the first
equation of (52), we get
$$
\frac{d^2r}{d\theta ^2}(\theta _i)=\frac 1{A_i^{n-2}}((n+1)A_i^{n-1}-K).
$$
Also $A_1\leq (\frac K{n+1})^{1/n+1}$ and $A_2\geq (\frac K{n+1})^{1/(n-1)}$
by (52). Thus $r\equiv (\frac K{n+1})^{1/(n-1)}$.\qed
\enddemo

We assume now that $K>0$ and consider the following initial values 
$$
\aligned
r(0)&=a>0, \\ 
r^{\prime }(0)&=0.
\endaligned \tag{53}
$$
So we get
$$
\lambda =a^{n+1}-\frac K2a^2.  \tag{54}
$$
The condition $\lambda \neq 0$ is equivalent to the following condition 
$$
a\neq (\frac K2)^{\frac 1{n-1}}.
$$
Set $f(x)=x^{n+1}-(\lambda +\frac K2x^2)$ for $x\in (0,\infty )$. It
is easy to see that $f$ is increasing strictly on $(a,\infty )$, provided
that 
$$
a>(\frac K{n+1})^{\frac 1{n-1}}.  \tag{55}
$$
Under the condition (55), we get 
$$
\frac{d^2r}{d\theta ^2}(0)=\frac 1{a^{n-2}}((n+1)a^{n-1}-K)>0.
$$
It follows that $\frac{dr}{d\theta }>0$ for $\theta \in (0,\varepsilon )$,
if $\varepsilon $ is small enough. So we have from (51) that 
$$
\frac{dr}{d\theta }=\frac{r\sqrt{r^{2(n+1)}-(\lambda +\frac K2r^2)^2}}{
\lambda +\frac K2r^2}\text{ },  \tag{56}
$$
for $\theta \in (0,\varepsilon )$. Since $f$ is increasing, it is easy to
see from (55) that $\frac{dr}{d\theta }>0$ for $\theta >0$. Let $\theta
_{\max }$ be the maximal value such that the solution exits on $[0,\theta
_{\max })$. From (56) we get 
$$
\theta _{\max }=\int_a^{r(\theta _{\max })}\frac{\lambda +\frac K2\rho ^2}{
\rho \sqrt{\rho ^{2(n+1)}-(\lambda +\frac K2\rho ^2)^2}}d\rho .  \tag{57}
$$
When $n>1$, the integral (57) converges. So $\theta _{\max }$ is finite and 
$r$ is an increasing function on $[0,\theta _{\max })$. By ODE theory, there
are two possibilities: (i) $\lim_{\theta \rightarrow \theta _{\max
}}r(\theta )=+\infty $ ; (ii) $\lim_{\theta \rightarrow \theta _{\max
}}r(\theta )=r_0<+\infty $, but $\lim_{\theta \rightarrow \theta _{\max
}}r^{\prime }(\theta )=\infty $.

We assert that the case (ii) will not occur. For this case we will have by
(48) and (54) that 
$$
a^{n+1}-\frac K2a^2+\frac K2r^2(\theta _{\max })=0.
$$
But this is impossible, because $r^2(\theta _{\max })>a^2$. Hence we only
have 
$$
\lim_{\theta \rightarrow \theta _{\max }}r(\theta )=+\infty .
$$
We note that the solution of (44) is invariant under the reflection 
$$
\theta \rightarrow -\theta .
$$
By the reflection, we can get a positive solution $r$ on $(-\theta _{\max
},\theta _{\max })$. Two lines $\theta =\pm \theta _{\max }$ are asymptotic
lines of the solution curve. It is easy to see that the solution curve has
infinite length with respect to the metric $\widetilde{g}$ given by (43). In
conclusion, we have

\proclaim {Theorem 4.3}\newline
(i) The only closed $H-$minimal Lagrangian submanifold invariant under the
action (42) is the Lagrangian submanifold corresponding $r\equiv const.$,
which is given by 
$$
\align
S^1\times S^n &\longrightarrow C^{n+1} \\ 
(e^{i\theta },x_1,...,x_{n+1})& \longmapsto (\frac K{n+1})^{1/(n-1)}e^{i\theta
}(x_1,....,x_{n+1}),
\endalign
$$
where $x_1^2+\cdots +x_{n+1}^2=1$.\newline
(ii)There are infinitely many non-trivial complete $H-$minimal Lagrangian
immersions of $R\times S^n$ into $C^{n+1}$ which are invariant under the
action (42).
\endproclaim
\vskip 0.3 true cm
\subhead {4.3. $T^n-$invariant $H-$minimal Lagrangian submanifolds}
\endsubhead
\vskip 0.3 true cm
Let $G\cong T^n$ be the group of diagonal matrices in $SU(n+1)$, so that
each $\gamma \in G$ acts on $C^{n+1}$ ($n>1$) by 
$$
\gamma :(z_1,z_2,...,z_{n+1})\mapsto (e^{i\theta _1}z_1,e^{i\theta
_2}z_2,...,e^{i\theta _{n+1}}z_{n+1})  \tag{58}
$$
for some $\theta _1,...,\theta _{n+1}\in R$ with $\theta _1+\cdots +\theta
_{n+1}=0$. The moment map of $G$ is 
$$
(z_1,z_2,...,z_{n+1})\mapsto -\frac
i2(|z_1|^2-|z_{n+1}|^2,|z_2|^2-|z_{n+1}|^2,...,|z_n|^2-|z_{n+1}|^2).
$$
As $G$ is abelian, $Z(g^{*})=g^{*}$. Let $c_i\in R$ such that $c_1\cdots
c_n\neq 0$ and set $c=-\frac i2(c_1,c_2,...,c_n)$. The level set $\mu
^{-1}(c)$ is given by 
$$
\mu
^{-1}(c)=
\{(z_1,...,z_{n+1}):|z_1|^2-|z_{n+1}|^2=c_1,...,|z_n|^2-|z_{n+1}|^2=c_n\}.
$$
So we can introduce the following parametrization of $\mu ^{-1}(c)/G:$
$$
\align
(r,\theta )\overset{F}\to {\longmapsto }(\sqrt{r^2+c_1},&\sqrt{r^2+c_2},...,\sqrt{
r^2+c_n},r\cos \theta +\sqrt{-1}r\sin \theta ), \\ 
&\sqrt{\sigma }\leq r<\infty ,\ 0\leq \theta <2\pi ,
\endalign
$$
where $\sigma =\max_{1\leq i\leq n}\{-c_i,0\}$. Now we hope to derive the
metric on $\mu ^{-1}(c)/G$. Let $\xi _i,i=1,...,n,$ be the standard basis of
the Lie algebra of $T^n$. Then the tangent space of the orbit at $F(r,\theta
)$ is spanned by 
$$
\aligned
\phi (\xi _1)=(iz_1,0,...,0,-iz_{n+1})&\longleftrightarrow (0,\sqrt{r^2+c_1}
,...,r\sin \theta ,-r\cos \theta ) \\ 
&....  \\ 
\phi (\xi _n)=(0,0,...,iz_n,-iz_{n+1})&\longleftrightarrow (0,0,...,0,\sqrt{
r^2+c_n},r\sin \theta ,-r\cos \theta )
\endaligned \tag{59} 
$$
$$
\aligned
dF(\frac \partial {\partial r})=&(\frac r{\sqrt{r^2+c_1}},0,\frac r{\sqrt{
r^2+c_2}},0,...,\frac r{\sqrt{r^2+c_n}},0,\cos \theta ,\sin \theta ), \\ 
dF(\frac \partial {\partial \theta })=&(0,0,...0,0,-r\sin \theta ,r\cos
\theta )\endaligned \tag{60}
$$
Set $\Phi _{ij}=\langle \phi (\xi _i),\phi (\xi _i)\rangle$, $i,j=1,...,n$. From (59), we
get 
$$
\aligned
\Phi _{ii}&=2r^2+c_i, \\ 
\Phi _{ij}&=r^2\qquad \text{for }i\neq j.
\endaligned \tag{61}
$$
The induce metric on each orbit is given by 
$$
ds_{T^n}^2=\Phi _{ij}d\theta ^id\theta ^j.  \tag{62}
$$
The volume of the orbit corresponding to the point $(r,\theta )$ is (up to a
constant) 
$$
V(r,\theta )=\sqrt{\det (\Phi _{ij})}.  \tag{63}
$$
By a direct computation, we have 
$$
\det \Phi _{ij}=\sum_{j=1}^nr^2(r^2+c_1)\cdots \widehat{(r^2+c_j)}\cdots
(r^2+c_n)+\prod_{k=1}^n(r^2+c_k).  \tag{64}
$$
Obviously $dF(\frac \partial {\partial r})\bot span\{\phi (\xi _1),...,\phi
(\xi _n)\}$ and 
$$
\langle dF(\frac \partial {\partial r}),dF(\frac \partial {\partial r})\rangle =\frac{\det
(\Phi _{ij})}{\prod_{j=1}^n(r^2+c_j)}.  \tag{65}
$$
Note that $dF(\frac \partial {\partial \theta })$ is not horizontal w.r.t. $
\pi :\mu ^{-1}(c)\rightarrow \mu ^{-1}(c)/G$. Denote by $Pj_H$ the
projection on the horizontal space of the fibration $\pi $. By an elementary
computation, we may get 
$$
|Pj_HdF(\frac \partial {\partial \theta })|^2=\frac{r^2\prod_{k=1}^n(r^2+c_k)
}{\det (\Phi _{ij})}.  \tag{66}
$$
Also $\langle dF(\frac \partial {\partial r}),dF(\frac \partial {\partial \theta
})\rangle =0$. So we get from (65) and (66) the induced Kaehler metric on the
orbit space as follows 
$$
\widetilde{g}=\frac{\det (\Phi _{ij})}{\prod_{j=1}^n(r^2+c_j)}dr^2+\frac{
r^2\prod_{k=1}^n(r^2+c_k)}{\det (\Phi _{ij})}d\theta ^2.  \tag{67}
$$
Thus the Hsiang-Lawson metric on $\mu ^{-1}(c)/G$ is given by (up to a
constant)
$$
\widetilde{g}_{HL}=\frac{[\det (\Phi _{ij})]^2}{\prod_{k=1}^n(r^2+c_k)}
dr^2+r^2\prod_{k=1}^n(r^2+c_k)d\theta ^2.  \tag{68}
$$
The $H-$minimal equation is 
$$
V^2k_{\widetilde{L}}=K,  \tag{69}
$$
where $K$ is a constant. If there exists a point $p\in \widetilde{L}$ such
that $k_{\widetilde{L}}(p)=0$, then we see that $k_{\widetilde{L}}\equiv 0$.
This corresponds the special Lagrangian submanifold found by Harvey-Lawson
[HL](see Remark 4.1). The points with $\det (\Phi _{ij})=0$ correspond to
the singular points on $\mu ^{-1}(c)/G$.

Let $\widetilde{L}$ be given by $r=r(\theta )$ on $\mu ^{-1}(c)/G-\{\det
(\Phi _{ij})=0\}$. By a similar method as in previous sections, we may show
that (69) is the Euler-Lagrange equation of the following functional:
$$
J=\int_{\theta _1}^{\theta _2}[\sqrt{\frac{[\det (\Phi _{ij})]^2}{
\prod_{k=1}^n(r^2+c_k)}(r^{\prime })^2+r^2\prod_{k=1}^n(r^2+c_k)}-\frac
K2r^2]d\theta .  \tag{70}
$$
So the Hamiltonian for the Euler-Lagrange equation of the functional (70) is
given by 
$$
H(\theta ,r,p)=-\frac{r\prod_{k=1}^n(r^2+c_k)}{\det (\Phi _{ij})}\sqrt{\frac{
[\det (\Phi _{ij})]^2}{\prod_{k=1}^n(r^2+c_k)}-p^2}+\frac K2r^2,  \tag{71}
$$
where 
$$
p=L_{r^{\prime }}=\frac{[\det (\Phi _{ij})]^2r^{\prime }}{
\prod_{k=1}^n(r^2+c_k)\sqrt{\frac{[\det (\Phi _{ij})]^2}{
\prod_{k=1}^n(r^2+c_k)}(r^{\prime })^2+r^2\prod_{k=1}^n(r^2+c_k)}}.  \tag{72}
$$
Since $H(\theta ,r,p)$ doesn't depend on $\theta $ explicitly, it must be a
constant of the motion. Thus we get from (71) and (72) that 
$$
\frac{r^2\prod_{k=1}^n(r^2+c_k)}{\sqrt{\frac{[\det (\Phi _{ij})]^2}{
\prod_{k=1}^n(r^2+c_k)}(r^{\prime })^2+r^2\prod_{k=1}^n(r^2+c_k)}}=\lambda
+\frac K2r^2.  \tag{73}
$$
From (73), we see that there exists no point $\theta _0$ such that 
$$
\lim_{\theta \rightarrow \theta _0}[r^2(\theta )\prod_{k=1}^n(r^2(\theta
)+c_k)]=0
$$
as $\theta \rightarrow \theta _0$, provided that $\lambda \neq 0,-\frac
K2c_i $ ($i=1,...,n$).

To simplified the discussion, we assume now that $K>0$ and consider the
following initial values 
$$
\aligned
r(0)&=a>\sqrt{\sigma }, \\ 
r^{\prime }(0)&=0.
\endaligned \tag{74}
$$
From (73) and (74), we get 
$$
\lambda =a\sqrt{\prod_{k=1}^n(a^2+c_k)}-\frac K2a^2.  \tag{75}
$$
In following, we always choose initial values such that $\lambda \neq
0,-\frac K2c_i$.

If $r^{\prime \prime }(0)=0$, then the solution $r(\theta )\equiv a$ by the
uniqueness Theorem of ODE. From (73), we have 
$$
\frac{dr}{d\theta }=\pm \frac{r\prod_{k=1}^n(r^2+c_k)\sqrt{
r^2\prod_{k=1}^n(r^2+c_k)-(\lambda +\frac K2r^2)^2}}{\det (\Phi
_{ij})(\lambda +\frac K2r^2)}.  \tag{76}
$$
Without loss of generality, we assume that $r^{\prime \prime }(0)>0$. By
taking derivative of (76) and using (74), we see that this assumption is
equivalent to 
$$
\frac{\det (\Phi _{ij}(a))}{a\sqrt{\prod_{k=1}^n(a^2+c_k)}}>K.  \tag{77}
$$
Set $F(x)=\sqrt{x\prod_{k=1}^n(x+c_k)}-(\lambda +\frac K2x)$ on $(\sigma,+\infty )
$. Obviously, we may choose $a$ large enough so that (77) is satisfied and $F
$ is increasing strictly on $(a^2,+\infty )$. Let $\theta _{\max }$ be the
maximal value such that the solution $r(\theta )$ exists on $[0,\theta
_{\max })$.  From (76), we get 
$$
\theta _{\max }=\int_a^{r(\theta _{\max })}\frac{\det (\Phi _{ij})(\lambda
+\frac K2r^2)dr}{r\prod_{k=1}^n(r^2+c_k)\sqrt{r^2\prod_{k=1}^n(r^2+c_k)-(
\lambda +\frac K2r^2)^2}}.  \tag{78}
$$
Since the integral (79) exists, $\theta _{\max }$ is finite. Similar to the
previous discussion in section 4.2, we see that $r(\theta )$ is a strictly
increasing solution on $[0,\theta _{\max })$ such that 
$$
\lim_{\theta \rightarrow \theta _{\max }}r(\theta )=+\infty .  \tag{79}
$$
Note also that the solution of (69) (or see (73)) is invariant under the
reflection 
$$
\theta \rightarrow -\theta .
$$
By the reflection, we can get a positive solution $r$ on $(-\theta _{\max
},\theta _{\max })$ so that 
$$
\lim_{\theta \rightarrow \pm \theta _{\max}}r(\theta )=+\infty .
$$

Now we show that the solution curves $r=r(\theta )$ are complete with
respect to the metric $\widetilde{g}$ given by (67). For the solution $r=r(\theta )$ on 
$(-\theta _{\max },\theta _{\max })$, we have from (67) and (64) that 
$$
\align
L_{\widetilde{g}}(r) &=2\int_0^{\theta _{\max }}\sqrt{\frac{\det (\Phi
_{ij})}{\prod_{j=1}^n(r^2+c_j)}(r^{\prime })^2+\frac{r^2
\prod_{k=1}^n(r^2+c_k)}{\det (\Phi _{ij})}}d\theta \\
&\geq 2\int_0^{\theta _{\max }}\sqrt{\frac{\det (\Phi _{ij})}{
\prod_{j=1}^n(r^2+c_j)}(r^{\prime })^2}d\theta \\
&\geq 2\int_0^{\theta _{\max }}r^{\prime }d\theta \\
&=2(r(\theta _{\max })-a) \\
&=+\infty .
\endalign
$$
In conclusion, we have proved that

\proclaim {Theorem 4.4}
There are infinitely many non-trivial complete $H-$minimal immersions of $
R\times T^n$ into $C^{n+1}$, which are invariant under the action (58).
\endproclaim

\remark {Remark 4.1}\roster
\item Note that the $H-$minimal immersions of $R\times T^n$ in Theorem 4.4 are different from
those in Theorem 4.1, because their actions of groups are different; 
\item If we set $R=\sqrt{r^2\prod_{j=1}^n(r^2+c_j)}$, then the Hsiang-Lawson
metric becomes
$$
\widetilde{g}_{HL}=dR^2+R^2d\theta ^2
$$
which is actually a flat metric. Using the parametrization here, the $T^n-$
invariant special Lagrangian $(n+1)-$folds in $C^{n+1}$ of Harvey-Lawson
[HL] correspond to the straight lines $R\sin \theta \equiv const.$ or $R\cos
\theta \equiv const.$ on the orbit space $\mu ^{-1}(c)/G$. Up to a $SO(2)-$
motion, the different straight lines in the $(x=R\cos \theta$, $y=R\sin \theta)-$plane 
corresponds to special Lagrangian $(n+1)-$folds in $C^{n+1}$
with different phases. So we see that the complete $H-$minimal submanifolds
in Theorem 4.4 are asymptotic to two singular $T^n-$invariant $\pm\theta_0-$special
Lagrangian $(n+1)-$folds.
\endroster
\endremark
\vskip 0.3 true cm
{\bf Acknowledgments:}\quad Partial work of this paper was done while the
author visited Math. Dept. of Harvard University. The author would like to
thank Prof. S.T. Yau for his enthusiastic help and valuable suggestions. He
would also like to thank Professors C.H. Gu and H.S. Hu for their constant
encouragement and helpful comments; thank Prof. G.F. Wang, Dr. X.W. Wang and
Dr. Q. C. Ji for their valuable conversations. \footnote{
The author announced the paper on the International conference on Geometric analysis  at Beijing in 2003. 
The method of this paper was applied by Dr. Q.C. Ji to the complex hyperbolic space. Before submitting
the paper to Math. Arxiv, the author noticed a recent different
constructions of H-minimal submanifolds by Castro-Li-Urbanno(DG/0412046).}

\vskip 1 true cm

\vskip 1 true cm

Institute of Mathematics

Fudan University

Shanghai, 200433

P.R.China
\vskip 0.2 true cm
yxdong\@fudan.edu.cn

\enddocument